\documentclass[11pt]{article}

\usepackage{latexsym}
\usepackage{amssymb}

\newtheorem{thm}{Theorem}
\newtheorem{dfn}{Definition}

\begin{document}
{
\begin{center}
{\Large\bf
On the density of polynomials in some $L^2(M)$ spaces.
}
\end{center}
\begin{center}
{\bf S.M. Zagorodnyuk}
\end{center}

\section{Introduction.}
In this paper we shall study the density of polynomials in some $L^2(M)$ spaces. Two choices of
the measure $M$ and polynomials will be considered:
\begin{itemize}
\item[(A)]
a $\mathbb{C}_{N\times N}^\geq$-valued measure $M$ on $\mathfrak{B}(\mathbb{R})$ and vector-valued polynomials:
\begin{equation}
\label{f1_1}
p(x) = (p_0(x),p_1(x),...,p_{N-1}(x)),
\end{equation}
where $p_j(x)$ are complex polynomials, $0\leq j\leq N-1$; $N\in \mathbb{N}$;
\item[(B)]
a scalar non-negative Borel measure $\sigma$ in a strip
\begin{equation}
\label{f1_2}
\Pi = \{ (x,\varphi):\ x\in \mathbb{R},\ \varphi\in [-\pi,\pi) \},
\end{equation}
and power-trigonometric polynomials:
\begin{equation}
\label{f1_3}
p(x,\varphi) = \sum_{m=0}^\infty \sum_{n=-\infty}^\infty \alpha_{m,n} x^m e^{in\varphi},\ \alpha_{m,n}\in \mathbb{C},
\end{equation}
where all but finite number of coefficients $\alpha_{m,n}$ are zeros.
\end{itemize}
The case $(A)$ is closely related to the matrix Hamburger moment problem which consists of
finding a left-continuous non-decreasing matrix function $M(x) = ( m_{k,l}(x) )_{k,l=0}^{N-1}$
on $\mathbb{R}$, $M(-\infty)=0$, such that
\begin{equation}
\label{f1_4}
\int_\mathbb{R} x^n dM(x) = S_n,\qquad n\in \mathbb{Z}_+,
\end{equation}
where $\{ S_n \}_{n=0}^\infty$ is a prescribed sequence of Hermitian $(N\times N)$ complex matrices, $N\in \mathbb{N}$.
In the scalar case ($N=1$) it is well known that polynomials are dense in $L^2(M)$ on the real line
if and only if $M$ is a canonical solution of the corresponding moment problem~\cite{cit_900_Akh}.

\noindent
In the case of an arbitrary $N$ and if the matrix Hamburger moment problem
is completely indetermined, the density of polynomials is equivalent to the fact that $M$
is a canonical solution of the moment problem~(\ref{f1_4}) (i.e. it corresponds to a
constant unitary matrix in the Nevanlinna type parameterization for solutions of~(\ref{f1_4}))~\cite{cit_1100_L}.

\noindent
On the other hand, the case (B) is related to the Devinatz moment problem:
to find a non-negative Borel measure $\mu$ in a strip $\Pi$
such that
\begin{equation}
\label{f1_5}
\int_\Pi x^m e^{in\varphi} d\mu = s_{m,n},\qquad m\in \mathbb{Z}_+, n\in \mathbb{Z},
\end{equation}
where $\{ s_{m,n} \}_{m\in \mathbb{Z}_+, n\in \mathbb{Z}}$ is a prescribed sequence of complex
numbers~\cite{cit_1150_Z}.

\noindent
In the both cases, we shall prove
that polynomials are dense in $L^2(M)$ if and only if $M$ is a canonical solution of the
corresponding moment problem, without any additional assumptions (definitions of the canonical
solutions shall be given below).
For this purpose, we derive a model for a finite set of commuting self-adjoint and unitary operators with
a spectrum of a finite multiplicity (precise definitions shall be stated below).
The latter is a generalization of the canonical model for a self-adjoint operator with
a spectrum of a finite multiplicity~\cite{cit_1200_AG}.
Using known descriptions of canonical solutions, we shall obtain conditions for the density of polynomials in $L^2(M)$.

{\bf Notations. } As usual, we denote by $\mathbb{R},\mathbb{C},\mathbb{N},\mathbb{Z},\mathbb{Z}_+$
the sets of real numbers, complex numbers, positive integers, integers and non-negative integers,
respectively; $\mathbb{C}_+ := \{ z\in \mathbb{C}:\ \frac{1}{2i} (z-\overline{z}) \geq 0 \}$.
By $\mathbb{C}_{n\times n}$ we denote a set of all $(n\times n)$ matrices with complex elements;
$\mathbb{C}_{n} := \mathbb{C}_{1\times n}$,\ $n\in \mathbb{N}$.
By $\mathbb{C}_{n\times n}^{\geq}$ we mean a set of all nonnegative Hermitian matrices from
$\mathbb{C}_{n\times n}$, $n\in \mathbb{N}$.
By $\mathbb{P}$ we denote a set of all complex polynomials. By $\mathbb{P}^N$ we mean a set
of vector-valued polynomials: $p(z) = (p_0(z),p_1(z),...,p_{N-1}(z))$; $p_j\in \mathbb{P}$,
$0\leq j\leq N-1$; $N\in \mathbb{N}$.
For a subset $S$ of the complex plane we denote by $\mathfrak{B}(S)$ the set of all Borel subsets of $S$.
Everywhere in this paper, all Hilbert spaces are assumed to be separable.
By $(\cdot,\cdot)_H$ and $\| \cdot \|_H$ we denote the scalar product and the norm in a Hilbert space $H$,
respectively. The indices may be omitted in obvious cases.
For a set $M$ in $H$, by $\overline{M}$ we mean the closure of $M$ in the norm $\| \cdot \|_H$. For
$\{ x_k \}_{k\in S}$, $x_k\in H$, we write
$\mathop{\rm Lin}\nolimits \{ x_k \}_{k\in S}$ for a set of linear combinations of vectors $\{ x_k \}_{k\in S}$
and $\mathop{\rm span}\nolimits \{ x_k \}_{k\in S} =
\overline{ \mathop{\rm Lin}\nolimits \{ x_k \}_{k\in S} }$.
Here $S$ is an arbitrary set of indices.
The identity operator in $H$ is denoted by $E=E_H$. For an arbitrary linear operator $A$ in $H$,
the operators $A^*$,$\overline{A}$,$A^{-1}$ mean its adjoint operator, its closure and its inverse
(if they exist). By $D(A)$ and $R(A)$ we mean the domain and the range of the operator $A$.
We denote by $R_z (A)$ the resolvent function of $A$, where $z$ belongs to the resolvent set of $A$.
If $A$ is bounded, then the norm of $A$ is denoted by $\| A \|$.
If $A$ is symmetric, we denote $\Delta_A(z) := (A - zE_H) D(A)$, $z\in \mathbb{C}$; and
$N_\lambda = N_\lambda(A) = H\ominus \Delta_A(\lambda)$, $\lambda\in \mathbb{C}\backslash \mathbb{R}$.
By $P^H_{H_1} = P_{H_1}$ we mean the operator of orthogonal projection in $H$ on a subspace
$H_1$ in $H$.

\noindent
We denote $D_{r,l} = \mathbb{R}^{r} \times [-\pi,\pi)^{l}
= \{ (x_1,x_2,...,x_r,\varphi_1,\varphi_2,...,\varphi_l),$
$x_j\in \mathbb{R},\ \varphi_k\in [-\pi,\pi),\ 1\leq j\leq r,\ 1\leq k\leq l \}$, $r,l\in \mathbb{Z}_+$.
Elements $u\in D_{r,l}$ we briefly denote by $u = (x,\varphi)$, $x=(x_1,x_2,...,x_r)$,
$\varphi = (\varphi_1,\varphi_2,...,\varphi_l)$. We mean $D_{r,0} = \mathbb{R}^{r}$;
$D_{0,l} = [-\pi,\pi)^{l}$.

\noindent
Let $M(\delta) = (m_{i,j}(\delta))_{i,j=0}^{N-1}$ be a
$\mathbb{C}_{N\times N}^\geq$-valued measure on $\mathfrak{B}(D_{r,l})$, and
$\tau = \tau_M (\delta) := \sum_{k=0}^{N-1} m_{k,k} (\delta)$;
$M'_\tau = (m'_{k,l})_{k,l=0}^{N-1} = ( dm_{k,l}/ d\tau_M )_{k,l=0}^{N-1}$; $N\in \mathbb{N}$.
We denote by $L^2(M)$ a set (of classes of equivalence)
of vector-valued functions
$f: D_{r,l}\rightarrow \mathbb{C}_N$, $f = (f_0,f_1,\ldots,f_{N-1})$, such that
(see, e.g.,~\cite{cit_1400_R},\cite{cit_10000_MM})
$$ \| f \|^2_{L^2(M)} := \int_{D_{r,l}}  f(u) \Psi(u) f^*(u) d\tau_M  < \infty. $$
The space $L^2(M)$ is a Hilbert space with the scalar product
$$ ( f,g )_{L^2(M)} := \int_{D_{r,l}}  f(u) \Psi(u) g^*(u) d\tau_M,\qquad f,g\in L^2(M). $$
Set
$$ W_n f(x,\varphi) = e^{i\varphi_n} f(x,\varphi),\qquad f\in L^2(M);\ 1\leq n\leq l; $$
and
$$ X_m f(x,\varphi) = x_m f(x,\varphi), $$
$$ f(x,\varphi)\in L^2(M):\ x_m f(x,\varphi)\in L^2(M);\
1\leq m\leq r. $$
Operators $W_n$ are unitary. In the usual manner~\cite{cit_12000_BS}, one can check that
operators $X_m$ are self-adjoint.

\section{A set of commuting self-adjoint and unitary operators with a spectrum of a finite
multiplicity.}
It is well known that a self-adjoint operator with a spectrum of a finite multiplicity
in a Hilbert space $H$ has a canonical model as a multiplication by an independent variable in $L^2(M)$.
Here $M$ is a $\mathbb{C}_{N\times N}^\geq$-valued measure on $\mathfrak{B}(\mathbb{R})$,
and $N$ is the multiplicity of the spectrum of $A$ \cite{cit_1200_AG}.
For our investigation on the density of polynomials, mentioned in the Introduction, we shall use
a generalization of this result to the case of an arbitrary finite set of commuting
self-adjoint and unitary operators. Moreover, we shall need a result which is a little more
general even in the classical case. Our method of proof is little different from the classical one
(we shall not use Lemma in~\cite[p.287]{cit_1200_AG}).

Consider a set
\begin{equation}
\label{f2_1}
\mathcal{A} = (S_1,S_2,...,S_\mathbf{r},U_1,U_2,...,U_\mathbf{l}),\quad \mathbf{r},\mathbf{l}\in \mathbb{Z}_+:\
\mathbf{r}+\mathbf{l}\not=0,
\end{equation}
where $S_j$ are self-adjoint operators and $U_k$ are unitary operators in a Hilbert space $H$,
$1\leq j\leq \mathbf{r}$, $1\leq k\leq \mathbf{l}$.
In the case $\mathbf{r}=0$ operators $S_j$ disappear. Analogously, for $\mathbf{l}=0$ we only have operators
$S_j$.
The set $\mathcal{A}$ is said to be a {\bf $SU$-set of order $(\mathbf{r},\mathbf{l})$}.

\noindent
The set $\mathcal{A}$ is called {\bf commuting} if operators $S_j$,$U_k$ pairwise commute. This mean that
\begin{equation}
\label{f2_2}
U_k U_m = U_m U_k,\quad 1\leq k,m\leq \mathbf{l};
\end{equation}
\begin{equation}
\label{f2_3}
U_k S_j \subset S_j U_k,\quad 1\leq j\leq \mathbf{r};\ 1\leq k\leq \mathbf{l};
\end{equation}
and the spectral measures of $S_j$ pairwise commute~\cite{cit_12000_BS}.
In this case, there exists a spectral measure $E(\delta)$, $\delta\in \mathfrak{B}( D_{\mathbf{r},\mathbf{l}} )$,
such that~\cite{cit_12000_BS}:
\begin{equation}
\label{f2_3_1}
S_j = \int_{ D_{\mathbf{r},\mathbf{l}} } x_j dE,\quad 1\leq j\leq \mathbf{r};
\end{equation}
\begin{equation}
\label{f2_3_2}
U_k = \int_{ D_{\mathbf{r},\mathbf{l}} } e^{i\varphi_k} dE,\quad 1\leq k\leq \mathbf{l}.
\end{equation}
We shall call $E$ {\bf the spectral measure of the commuting $SU$-set $\mathcal{A}$ of
order $(\mathbf{r},\mathbf{l})$}.

\noindent
We shall say that a commuting $SU$-set $\mathcal{A}$ of order $(\mathbf{r},\mathbf{l})$
{\bf has a spectrum of multiplicity $d$}, if
\begin{itemize}
\item[1)] there exist vectors $h_0,h_1,...,h_{d-1}$ in $H$ such that
\begin{equation}
\label{f2_4}
h_i \in D(S_1^{m_1}S_2^{m_2}...S_\mathbf{r}^{m_\mathbf{r}}),\quad m_1,m_2,...,m_\mathbf{r}\in \mathbb{Z}_+,\
0\leq i\leq d-1;
\end{equation}
$$ \mathop{\rm span}\nolimits \{
U_1^{n_1} U_2^{n_2}...U_l^{n_\mathbf{l}} S_1^{m_1}S_2^{m_2}...S_\mathbf{r}^{m_\mathbf{r}} h_i, $$
\begin{equation}
\label{f2_5}
m_1,m_2,...,m_\mathbf{r}\in \mathbb{Z}_+;\ n_1,n_2,...,n_\mathbf{r}\in \mathbb{Z};\
0\leq i\leq d-1 \} = H;
\end{equation}

\item[2)] ({\it minimality})
For arbitrary $\widetilde d\in \mathbb{Z}_+:\ \widetilde d < d$, and arbitrary
$\widetilde h_0,\widetilde h_1,...,\widetilde h_{d-1}$ in $H$,
at least one of conditions~(\ref{f2_4}),(\ref{f2_5}), with $\widetilde d$ instead of $d$, and
$\widetilde h_i$ instead of $h_i$, is not satisfied.

\end{itemize}
In the case $\mathbf{r}=0$, condition~(\ref{f2_4}) is redundant. Condition~(\ref{f2_5}) in
cases $\mathbf{r}=0$, $\mathbf{l}=0$, has no $U_k$ or $S_j$, respectively.

Set
$$ \vec e_i = (\delta_{0,i},\delta_{1,i},...,\delta_{N-1,i}),\qquad 0\leq i\leq N-1. $$

\begin{thm}
\label{t2_1}
Let $\mathcal{A}$ be a commuting $SU$-set of order $(\mathbf{r},\mathbf{l})$ in a Hilbert space $H$
which has a spectrum of multiplicity $d$.
Let $x_0,x_1,...,x_{N-1}$, $N\geq d$, be elements of $H$ such that
\begin{equation}
\label{f2_9}
x_i \in D(S_1^{m_1}S_2^{m_2}...S_\mathbf{r}^{m_\mathbf{r}}),\quad m_1,m_2,...,m_\mathbf{r}\in \mathbb{Z}_+,\
0\leq i\leq N-1;
\end{equation}
$$ \mathop{\rm span}\nolimits \{
U_1^{n_1} U_2^{n_2}...U_l^{n_\mathbf{l}} S_1^{m_1}S_2^{m_2}...S_\mathbf{r}^{m_\mathbf{r}} x_i, $$
\begin{equation}
\label{f2_10}
m_1,m_2,...,m_\mathbf{r}\in \mathbb{Z}_+;\ n_1,n_2,...,n_\mathbf{r}\in \mathbb{Z};\
0\leq i\leq N-1 \} = H.
\end{equation}
Set
\begin{equation}
\label{f2_11}
M(\delta) = \left(
(E(\delta) x_i,x_j)_H
\right)_{i,j=0}^{N-1},\qquad \delta\in \mathfrak{B}( D_{\mathbf{r},\mathbf{l}} ),
\end{equation}
where $E$ is the spectral measure of $\mathcal{A}$.

\noindent
Then there exists a unitary transformation $V$ which maps $L^2(M)$ onto $H$ such that:
\begin{equation}
\label{f2_12}
V^{-1} S_j V = X_j,\qquad 1\leq j\leq \mathbf{r};
\end{equation}
\begin{equation}
\label{f2_13}
V^{-1} U_k V = W_k,\qquad 1\leq k\leq \mathbf{l}.
\end{equation}
Moreover, we have
\begin{equation}
\label{f2_14}
V \vec e_s = x_s,\qquad 0\leq s\leq N-1.
\end{equation}
\end{thm}
{\bf Remark.} In the case $\mathbf{r}=0$ relations~(\ref{f2_9}),(\ref{f2_12}) should be removed, and
in~(\ref{f2_10}) operators $S_j$ disappear. In the case $\mathbf{l}=0$ relation~(\ref{f2_13})
should be removed and in~(\ref{f2_10}) operators $U_k$ disappear.

\noindent
{\bf Proof.}
Let $\chi_\delta(u)$ be the characteristic function of a set $\delta\in \mathfrak{B}( D_{\mathbf{r},\mathbf{l}} )$.
In the space $L^2(M)$ consider the following set:
\begin{equation}
\label{f2_15}
L := \mathop{\rm Lin}\nolimits \{ \chi_\delta(u) \vec e_s,\ \delta \in \mathfrak{B}( D_{\mathbf{r},\mathbf{l}} ),\
0\leq s\leq N-1 \}.
\end{equation}
Choose two arbitrary functions
\begin{equation}
\label{f2_16}
f(u) = \sum_{j=0}^{N-1} \sum_{\delta\in I_j} \alpha_j(\delta) \chi_\delta(u) \vec e_j,\quad
\alpha_j(\delta)\in \mathbb{C},
\end{equation}
\begin{equation}
\label{f2_17}
g(u) = \sum_{s=0}^{N-1} \sum_{\delta'\in J_s} \beta_s(\delta') \chi_{\delta'}(u) \vec e_s,\quad
\beta_s(\delta')\in \mathbb{C},
\end{equation}
where $I_j$,$J_s$ are some finite subsets of $\mathfrak{B}( D_{\mathbf{r},\mathbf{l}} )$.
We may write
$$ (f(u),g(u))_{L^2(M)} =
\sum_{j,s=0}^{N-1} \sum_{\delta\in I_j} \sum_{\delta'\in J_s} \alpha_j(\delta)
\overline{ \beta_s(\delta') } \int_{ D_{\mathbf{r},\mathbf{l}} } \chi_{\delta\cap\delta'} (u)
\vec e_j M'_\tau(u) \vec e_s^* d\tau_M $$
\begin{equation}
\label{f2_18}
\sum_{j,s=0}^{N-1} \sum_{\delta\in I_j} \sum_{\delta'\in J_s} \alpha_j(\delta)
\overline{ \beta_s(\delta') } m_{j,s} (\delta\cap\delta').
\end{equation}
Set
\begin{equation}
\label{f2_19}
x_f = \sum_{j=0}^{N-1} \sum_{\delta\in I_j} \alpha_j(\delta) E(\delta) x_j,\quad
x_g = \sum_{s=0}^{N-1} \sum_{\delta'\in J_s} \beta_s(\delta') E(\delta') x_s.
\end{equation}
Then
$$ (x_f,x_g)_H = \sum_{j,s=0}^{N-1} \sum_{\delta\in I_j} \sum_{\delta'\in J_s} \alpha_j(\delta)
\overline{ \beta_s(\delta') } (E(\delta) x_r, E(\delta') x_s)_H $$
\begin{equation}
\label{f2_21}
= \sum_{j,s=0}^{N-1} \sum_{\delta\in I_j} \sum_{\delta'\in J_s} \alpha_j(\delta)
\overline{ \beta_s(\delta') } m_{j,s} (\delta\cap\delta').
\end{equation}
Comparing relations~(\ref{f2_18}) and~(\ref{f2_21}) we obtain:
\begin{equation}
\label{f2_22}
(f,g)_{L^2(M)} = (x_f,x_g)_H.
\end{equation}
Now assume that $f$ and $g$ belong to the same class of equivalence in $L^2(M)$:
$\| f-g \|_{L^2(M)} = 0$. Then
$$ \| x_f - x_g \|_H^2 =
\left\|
\sum_{j=0}^{N-1} \left(
\sum_{\delta\in I_j} \alpha_j(\delta) E(\delta) -
\sum_{\delta\in J_j} \beta_j(\delta) E(\delta)
\right)
x_j
\right\|_H^2 $$
$$ =
\left\|
\sum_{j=0}^{N-1} \sum_{\delta\in I_j\cup J_j} c_j(\delta) E(\delta) x_j
\right\|_H^2, $$
where
\begin{equation}
\label{f2_23}
c_j(\delta) = \left\{\begin{array}{ccc} \alpha_j(\delta), & \delta\in I_j\backslash J_j\\
-\beta_j(\delta), & \delta\in J_j\backslash I_j\\
\alpha_j(\delta) - \beta_j(\delta), & \delta\in I_j\cap J_j\end{array}
\right..
\end{equation}
Set
\begin{equation}
\label{f2_24}
w(u) = \sum_{j=0}^{N-1} \sum_{\delta\in I_j\cup J_j} c_j(\delta) \chi_\delta(u) \vec e_j.
\end{equation}
Applying relation~(\ref{f2_22}) with $f=g=w$ we obtain:
$$ \| x_f - x_g \|_H^2 =
\| x_w \|_H^2 = \| w \|_{L^2(M)}^2 $$
$$ =
\left\|
\sum_{j=0}^{N-1}
\left(
\sum_{\delta\in I_j} \alpha_j(\delta) \chi_\delta(u) -
\sum_{\delta\in J_j} \beta_j(\delta) \chi_\delta(u)
\right) \vec e_j
\right\|_{L^2(M)}^2 =
\| f-g \| _{L^2(M)}^2 = 0. $$
Therefore a transformation $V$: $Vf = x_f$, is correctly defined on $L$, and $R(V)\subseteq H$.
Moreover, relation~(\ref{f2_22}) shows  that $V$ is an isometric transformation.
Since simple functions are dense in $L^2(M)$ (\cite[Theorem 3.11]{cit_1400_R}),
we have $\overline{L} = L^2(M)$. By continuity we extend $V$ on the whole $L^2(M)$.

\noindent
Suppose that $R(V)\not=H$. Then there exists $0\not= h\in H$, such that
$$ (E(\delta) x_s, h)_H = 0,\qquad \delta\in \mathfrak{B}(D_{\mathbf{r},\mathbf{l}}),\ 0\leq s\leq N-1. $$
Therefore we may write
$$ (U_1^{n_1} U_2^{n_2}...U_\mathbf{l}^{n_\mathbf{l}} S_1^{m_1}S_2^{m_2}...S_\mathbf{r}^{m_\mathbf{r}} x_s, h)_H $$
$$ = \int_{ D_{\mathbf{r},\mathbf{l}} } x_1^{m_1} x_2^{m_2} ... x_\mathbf{r}^{m_\mathbf{r}}
e^{ in_1\varphi_1 } e^{ in_2\varphi_2 } ... e^{ in_\mathbf{l}\varphi_l } d(E x_s, h)_H = 0, $$
$$ m_1,m_2,...,m_\mathbf{r}\in \mathbb{Z}_+,\ n_1,n_2,...,n_\mathbf{l}\in \mathbb{Z}. $$
By~(\ref{f2_10}) we get $h=0$. This contradiction proves that $R(V)=H$.
Thus, $V$ is a unitary transformation of $L^2(M)$ onto $H$.
Observe that relation~(\ref{f2_14}) holds.
Set
$$ L^2_i(M) = \{ f(u) = (f_0(u),f_1(u),...,f_{N-1}(u))\in L^2(M): $$
\begin{equation}
\label{f2_25}
\int_{ D_{\mathbf{r},\mathbf{l}} }
|f_s(u)|^2 dm_{s,s} < \infty,\ 0\leq s\leq N-1 \}.
\end{equation}
Here, as usual, we mean that $L^2_i(M)$ consists of classes of equivalence from $L^2(M)$, which have
at least one representative $f$ with square integrable components.
Observe that simple functions belong to $L^2_{i}(M)$ and therefore $L^2_{i}(M)$ is dense in $L^2(M)$.
Let us check that
\begin{equation}
\label{f2_26}
V f = \sum_{s=0}^{N-1} \int_{ D_{\mathbf{r},\mathbf{l}} }
f_s(u) dE x_s,\qquad f = (f_0,f_1,...,f_{N-1})\in L^2_i(M).
\end{equation}
Choose an arbitrary function $f = (f_0,f_1,...,f_{N-1})\in L^2_i(M)$. Let
\begin{equation}
\label{f2_27}
f_s^{k}(u) = \sum_{\delta\in I_{s,k}} \alpha_{s,k}(\delta) \chi_\delta(u),\quad 0\leq s\leq N-1;\ k\in \mathbb{N},
\end{equation}
where $I_{s,k}$ is a finite subset of $\mathfrak{B}(D_{\mathbf{r},\mathbf{l}})$,
be simple functions such that
\begin{equation}
\label{f2_28}
\int_{ D_{\mathbf{r},\mathbf{l}} }
|f_s(u) - f_s^k(u)|^2 dm_{s,s} \leq \frac{1}{k^2},\qquad 0\leq s\leq N-1;\ k\in \mathbb{N}.
\end{equation}
Then
\begin{equation}
\label{f2_29}
\| f(u) - \sum_{s=0}^{N-1} f_s^k(u) \vec e_s \|_{L^2(M)} \leq \frac{N}{k},\qquad  k\in \mathbb{N}.
\end{equation}
Set
$$ f^k(u) = \sum_{s=0}^{N-1} f_s^k(u) \vec e_s =
\sum_{s=0}^{N-1} \sum_{\delta\in I_{s,k}} \alpha_{s,k}(\delta) \chi_\delta(u) \vec e_s,\qquad
k\in \mathbb{N}. $$
Then
\begin{equation}
\label{f2_30}
\| f - f^k \|_{L^2(M)} \rightarrow 0,\quad \mbox{as } k\rightarrow\infty.
\end{equation}
Therefore
\begin{equation}
\label{f2_31}
\| Vf - Vf^k \|_{H} \rightarrow 0,\quad \mbox{as } k\rightarrow\infty.
\end{equation}
Observe that
\begin{equation}
\label{f2_32}
V f^k(u) = \sum_{s=0}^{N-1} \sum_{\delta\in I_{s,k}} \alpha_{s,k}(\delta) E(\delta) x_s,\qquad
k\in \mathbb{N}.
\end{equation}
We may write
$$ \left\|
\sum_{s=0}^{N-1} \int_{ D_{\mathbf{r},\mathbf{l}} }
f_s(u) dE x_s -
\sum_{s=0}^{N-1} \sum_{\delta\in I_{s,k}} \alpha_{s,k}(\delta) E(\delta) x_s
\right\|_H $$
$$ = \left\|
\sum_{s=0}^{N-1} \int_{ D_{\mathbf{r},\mathbf{l}} }
\left( f_s(u) - \sum_{\delta\in I_{s,k}} \alpha_{s,k}(\delta) \chi_\delta(u)) \right)
dE x_s
\right\|_H $$
$$ \leq
\sum_{s=0}^{N-1} \left\|
\int_{ D_{\mathbf{r},\mathbf{l}} }
\left( f_s(u) - \sum_{\delta\in I_{s,k}} \alpha_{s,k}(\delta) \chi_\delta(u) \right)
dE x_s
\right\|_H $$
$$ =
\sum_{s=0}^{N-1}
\left\{
\int_{ D_{\mathbf{r},\mathbf{l}} }
\left| f_s(u) - \sum_{\delta\in I_{s,k}} \alpha_{s,k}(\delta) \chi_\delta(u) \right|^2
d(Ex_s,x_s)_H
\right\}^{ \frac{1}{2} } \leq \frac{N}{k},\ k\in \mathbb{N}. $$
By the uniqueness of the limit we conclude that relation~(\ref{f2_26}) holds.

In the case $\mathbf{r}=0$, the following considerations until relations~(\ref{f2_38_1}),(\ref{f2_38_2})
are redundant, and in these relations one should choose $f\in L^2_{i}(M)$.

Set
$$ L^2_{i;2}(M) = \{ f(x,\varphi) = (f_0(x,\varphi),f_1(x,\varphi),...,f_{N-1}(x,\varphi))\in L^2(M): $$
$$ \int_{ D_{\mathbf{r},\mathbf{l}} }
|f_s(x,\varphi)|^2 dm_{s,s} < \infty,\
\int_{ D_{\mathbf{r},\mathbf{l}} }
|x_k f_s(x,\varphi)|^2 dm_{s,s} < \infty, $$
\begin{equation}
\label{f2_35}
1\leq k\leq \mathbf{r},\
0\leq s\leq N-1 \}.
\end{equation}
Of course, $L^2_{i;2}(M) \subseteq L^2_{i}(M)$, and $L^2_{i;2}(M) \subseteq D(X_k)$,
$1\leq k\leq \mathbf{r}$.
Moreover, we have
\begin{equation}
\label{f2_36}
X_m L^2_{i;2}(M) \subseteq L^2_{i}(M),\qquad 1\leq m\leq \mathbf{r}.
\end{equation}
Observe  that functions
\begin{equation}
\label{f2_37}
\chi_{\delta\cap\delta_k} (x,\varphi) \vec e_s,\quad \delta\in\mathfrak{B}(D_{\mathbf{r},\mathbf{l}}),\
0\leq s\leq N-1,
\end{equation}
\begin{equation}
\label{f2_38}
\delta_k = \{ (x,\varphi)\in D_{\mathbf{r},\mathbf{l}}:\ |x_m|\leq k,\ 1\leq m\leq \mathbf{r} \},\
k\in \mathbb{N},
\end{equation}
belong to $L^2_{i;2}(M)$. Therefore $L^2_{i;2}(M)$ is dense in $L^2(M)$.

Choose an arbitrary function $f\in L^2_{i;2}(M)$. By virtue of relation~(\ref{f2_26}) we may write:
\begin{equation}
\label{f2_38_1}
V f = \sum_{s=0}^{N-1} \int_{ D_{\mathbf{r},\mathbf{l}} } f_s(x,\varphi) dE x_s,
\end{equation}
$$ V X_m f = \sum_{s=0}^{N-1} \int_{ D_{\mathbf{r},\mathbf{l}} } x_m f_s(x,\varphi) dE x_s
= \sum_{s=0}^{N-1} S_m \int_{ D_{\mathbf{r},\mathbf{l}} } f_s(x,\varphi) dE x_s
= S_m V f, $$
\begin{equation}
\label{f2_38_2}
V W_n f = \sum_{s=0}^{N-1} \int_{ D_{\mathbf{r},\mathbf{l}} } e^{i\varphi_n} f_s(x,\varphi) dE x_s
= \sum_{s=0}^{N-1} U_n \int_{ D_{\mathbf{r},\mathbf{l}} } f_s(x,\varphi) dE x_s
= U_n V f,
\end{equation}
where $1\leq m\leq \mathbf{r}$, $1\leq n\leq \mathbf{l}$.
By continuity, from the latter relation we obtain that relation~(\ref{f2_13}) holds.
In the case $\mathbf{r}=0$ this completes the proof. In the opposite case
we may write
\begin{equation}
\label{f2_39}
X_m f = V^{-1} S_m V f,\qquad f\in L^2_{i;2}(M),\ 1\leq m\leq \mathbf{r}.
\end{equation}
Let us prove that
\begin{equation}
\label{f2_40}
L^2_{i;2}(M) \subseteq (X_m \pm i E_{ L^2(M) } ) L^2_{i;2}(M).
\end{equation}
Choose an arbitrary function $f=(f_0,f_1,...,f_{N-1})\in L^2_{i;2}(M)$.
Observe that
\begin{equation}
\label{f2_41}
g_\pm(x,\varphi) := \frac{1}{ x_m \pm i } (f_0(x,\varphi),f_1(x,\varphi),...,f_{N-1}(x,\varphi)) \in L^2_{i;2}(M).
\end{equation}
Therefore $(X_m \pm i E_{ L^2(M) } ) g_\pm(x,\varphi) = f$.
Thus, relation~(\ref{f2_40}) is true.
This relation means that operators $X_m$ and $V^{-1} S_m V$, restricted to $L^2_{i;2}(M)$,
are essentially self-adjoint. Therefore they have a unique self-adjoint extension. Since
operators $X_m$ and $V^{-1} S_m V$ are self-adjoint extensions, we conclude that relation~(\ref{f2_12}) holds.
$\Box$

\section{Density of polynomials: the case (A).}
Let $M = (m_{k,l})_{k,l=0}^{N-1}$ be a
$\mathbb{C}_{N\times N}^\geq$-valued measure on $\mathfrak{B}(\mathbb{R})$, $N\in \mathbb{N}$, such that
\begin{equation}
\label{f3_1}
\int_\mathbb{R} x^n dm_{k,l}\
\mbox{exist},\quad n\in \mathbb{Z}_+;\ 0\leq k,l\leq N-1.
\end{equation}
In this section, we shall use the same notation for matrix-valued measures $M(\delta)$ on $\mathfrak{B}(\mathbb{R})$
and their distribution functions $M(x)$, $x\in \mathbb{R}$~\cite{cit_10000_MM}.
Set
\begin{equation}
\label{f3_2}
S_n := \int_\mathbb{R} x^n dM,\qquad n\in \mathbb{Z}_+,
\end{equation}
and consider the matrix Hamburger moment problem with moments $\{ S_n \}_{n\in \mathbb{Z}_+}$.
Set
\begin{equation}
\label{f3_3}
\Gamma_n = (S_{k+l})_{k,l=0}^n,\ n\in \mathbb{Z}_+;\quad
\Gamma = (S_{k+l})_{k,l=0}^\infty = (\Gamma_{n,m})_{n,m=0}^\infty,\ \Gamma_{n,m}\in \mathbb{C}.
\end{equation}
Since the moment problem has a solution we have
$$ \Gamma_n \geq 0,\qquad n\in \mathbb{Z}_+. $$
There exists a Hilbert space $H$ and a sequence $\{ x_n \}_{n=0}^\infty$ in $H$,
such that $\mathop{\rm span}\nolimits \{ x_n \}_{n\in \mathbb{Z}_+} = H$, and~\cite{cit_14000_Z}
\begin{equation}
\label{f3_4}
(x_n,x_m)_H = \Gamma_{n,m},\qquad n,m\in \mathbb{Z}_+.
\end{equation}
Let $A$ be a linear operator with $D(A) = \mathop{\rm Lin}\nolimits \{ x_n \}_{n\in \mathbb{Z}_+}$, defined by
equalities
$$ A x_k =  x_{k+N},\qquad k\in \mathbb{Z}_+. $$
In~\cite{cit_14000_Z} it was shown that $A$ is a correctly defined symmetric operator in $H$.
Denote by $\mathbf{F} = \mathbf{F}(\overline{A})$ a set of all analytic in $\mathbb{C}_+$
operator-valued functions $F(\lambda)$, which values are contractions
which map $N_i(\overline{A})$ into $N_{-i}(\overline{A})$ ($\| F(\lambda) \|\leq 1$).
In~\cite[Theorem 4]{cit_14000_Z} it was proved that all solutions of the moment problem have the following form:
\begin{equation}
\label{f3_5}
\mathbf{M}(x) = (\mathbf{m}_{k,j} (x))_{k,j=0}^{N-1},
\end{equation}
where $\mathbf{m}_{k,j}$ satisfy the following relation
\begin{equation}
\label{f3_6}
\int_\mathbb{R} \frac{1}{x-\lambda} d \mathbf{m}_{k,j} (x) =
( (A_{F(\lambda)} - \lambda E_H)^{-1} x_k, x_j)_H,\qquad \lambda\in \mathbb{C}_+,
\end{equation}
where  $A_{F(\lambda)}$ is the
quasiself-adjoint extension of $\overline{A}$ defined by $F(\lambda)\in \mathbf{F}(\overline{A})$.

On the other hand, to any operator function $F(\lambda)\in \mathbf{F}(\overline{A})$  there corresponds by
relation~(\ref{f3_6}) a solution of the matrix Hamburger moment problem.
The correspondence between all operator functions $F(\lambda)\in \mathbf{F}(\overline{A})$ and all solutions
of the moment problem, established by relation~(\ref{f3_6}), is bijective.

Relation~(\ref{f3_6}) may be written in the following form:
\begin{equation}
\label{f3_7}
\int_\mathbb{R} \frac{1}{x-\lambda} d \mathbf{m}_{k,j} (x) =
( \mathbf{R}_\lambda (\overline{A}) x_k, x_j)_H,\qquad \lambda\in \mathbb{C}_+,
\end{equation}
where $\mathbf{R}_\lambda (\overline{A})$ is a generalized resolvent of $\overline{A}$.
The correspondence between all generalized resolvents and all solutions of the moment problem is
bijective.

\noindent
From relation~(\ref{f3_7}) it follows that~(\cite[Theorem 2]{cit_14000_Z})
\begin{equation}
\label{f3_8}
\mathbf{M}(t) = (\mathbf{m}_{k,j} (t) )_{k,j=0}^{N-1},\quad
\mathbf{m}_{k,j} (t) = ( \mathbf{E}_t x_k, x_j)_H,\qquad t\in \mathbb{R},
\end{equation}
where $\mathbf{E}_t$ is a spectral function of $\overline{A}$. The latter means that
$\mathbf{E}_t = P^{ \widehat H } _H \widehat E_t$, where $\widehat E_t$ is
the orthogonal resolution of unity of a self-adjoint operator $\widehat A\supseteq A$ in a Hilbert space
$\widehat H\supseteq H$.
The correspondence between all spectral functions and all solutions of the moment problem is
bijective, as well.
\begin{dfn}
A solution $\mathbf{M}(t) = (\mathbf{m}_{k,j} (t) )_{k,j=0}^{N-1}$ of the matrix Hamburger
moment problem~(\ref{f1_4}) is said to be {\bf canonical}, if it corresponds by relation~(\ref{f3_8})
to an orthogonal spectral function of $\overline{A}$, i.e. to a spectral function generated by a self-adjoint
extension $\widehat A\supseteq \overline{A}$ inside $H$.
\end{dfn}
From this definition we see that {\it canonical solutions exist if and only if the defect numbers of $A$ are equal}.
Observe that $\mathbf{M}(t) = (\mathbf{m}_{k,j} (t) )_{k,j=0}^{N-1}$ is a canonical solution of the matrix Hamburger
moment problem~(\ref{f1_4}) if and only if it corresponds to an orthogonal resolvent of $\overline{A}$,
i.e. to a usual resolvent of
a self-adjoint extension $\widehat A\supseteq \overline{A}$ inside $H$, in relation~(\ref{f3_7}).
Assume that the defect numbers of $A$ are equal.
From the Shtraus formula for generalized
resolvents~\cite[Theorem 7]{cit_25000_S}, it easily follows that
the orthogonal resolvents of $\overline{A}$ correspond to $F(\lambda)\equiv C$, $C$ is a unitary operator from
$N_i(\overline{A})$ onto $N_{-i}(\overline{A})$.
Consequently, canonical solutions of the moment problem correspond in relation~(\ref{f3_6})
to functions $F(\lambda)\equiv C$, $C$ is a unitary operator from
$N_i(\overline{A})$ onto $N_{-i}(\overline{A})$.

\begin{thm}
\label{t3_1}
Let $M = (m_{k,l})_{k,l=0}^{N-1}$ be a
$\mathbb{C}_{N\times N}^\geq$-valued measure on $\mathfrak{B}(\mathbb{R})$, $N\in \mathbb{N}$, such that
relation~(\ref{f3_1}) holds. Let $L^2_0(M)$ be the closure in $L^2(M)$ of a set of all
vector-valued polynomials $p\in \mathbb{P}^N$.
Consider the matrix Hamburger moment problem with moments $\{ S_n \}_{n\in \mathbb{Z}_+}$
defined by~(\ref{f3_2}). Consider
a Hilbert space $H$ and a sequence $\{ x_n \}_{n=0}^\infty$ in $H$,
such that $\mathop{\rm span}\nolimits \{ x_n \}_{n\in \mathbb{Z}_+} = H$, and
relation~(\ref{f3_4}) holds.
Let $A$ be a linear operator with $D(A) = \mathop{\rm Lin}\nolimits \{ x_n \}_{n\in \mathbb{Z}_+}$, defined by
equalities
$$ A x_k =  x_{k+N},\qquad k\in \mathbb{Z}_+. $$
The following conditions are equivalent:
\begin{itemize}
\item[(i)] $L^2_0(M) = L^2(M)$;

\item[(ii)] $M$ is a canonical solution of the corresponding matrix Hamburger moment problem;

\item[(iii)] $M(x) = (m_{k,j} (x))_{k,j=0}^{N-1}$ satisfy the following relation:
\begin{equation}
\label{f3_9}
\int_\mathbb{R} \frac{1}{x-\lambda} d m_{k,j} (x) =
( (A_{U} - \lambda E_H)^{-1} x_k, x_j)_H,\qquad \lambda\in \mathbb{C}_+,
\end{equation}
where  $A_{U}$ is a
quasiself-adjoint extension of $\overline{A}$ defined by a unitary operator $U$ from
$N_i(\overline{A})$ onto $N_{-i}(\overline{A})$.
The latter is equivalent to the fact that $A_U$ is a self-adjoint extension of $A$ inside $H$.

\item[(iv)] For every $\lambda\in \mathbb{C}_+$, there exists a linear bounded operator $D_\lambda$
in $H$ such that
\begin{equation}
\label{f3_10}
( D_\lambda x_{Nk+r}, x_{Nl+s} )_H = \int_\mathbb{R} \frac{ x^{k+l} }{ x-\lambda } dm_{r,s},\
0\leq r,s\leq N-1;\ k,l\in \mathbb{Z}_+,
\end{equation}
which is invertible and
\begin{equation}
\label{f3_11}
D_\lambda^{-1} + \lambda E_H \equiv A_U,
\end{equation}
where $A_U$ is a self-adjoint extension of $A$ inside $H$.
\end{itemize}
\end{thm}
{\bf Proof. }
(i)$\Rightarrow$(ii):
Repeating arguments from~\cite[pp.276-278]{cit_14000_Z} we construct a self-adjoint extension
$\widehat A$ of $A$, which acts in $H\oplus (L^2(M)\ominus L^2_0(M)) = H$, and
\begin{equation}
\label{f3_12}
m_{k,j}(t) = (\widehat E_t x_k,x_j)_H,
\end{equation}
where $\widehat E_t$
is a left-continuous resolution of unity of $\widehat A$.
Thus, $M$ is a canonical solution of the moment problem.

\noindent
(ii)$\Rightarrow$(i):
Let $M=(m_{k,j})_{k,j=0}^{N-1}$ has form~(\ref{f3_12}), where $\widehat E_t$
is a left-continuous resolution of unity of a self-adjoint operator $\widehat A\supseteq A$ in $H$.
Since $\widehat A x_n = A x_n = x_{n+N}$, $n\in \mathbb{Z}_+$, then by the induction argument we get
\begin{equation}
\label{f3_13}
\widehat A^r x_s = x_{rN+s},\qquad 0\leq s\leq N-1;\ r\in \mathbb{Z}_+.
\end{equation}
Therefore
\begin{equation}
\label{f3_14}
\mathop{\rm span}\nolimits \{ A^r x_s,\ 0\leq s\leq N-1;\ r\in \mathbb{Z}_+ \} = H.
\end{equation}
Thus, $\widehat A$ has a spectrum of multiplicity $d\leq N$.
By Theorem~\ref{t2_1} there exists a unitary transformation $W$ which maps $L^2(M)$ onto $H$ such that:
\begin{equation}
\label{f3_15}
W^{-1} \widehat A W = X,
\end{equation}
\begin{equation}
\label{f3_16}
W \vec e_s = x_s,\qquad 0\leq s\leq N-1,
\end{equation}
where $X$ is the operator of multiplication by an independent variable in $L^2(M)$.
Let us check that
\begin{equation}
\label{f3_17}
W x^k \vec e_s = x_{kN+s},\qquad 0\leq s\leq N-1;\ k\in \mathbb{Z}_+.
\end{equation}
Fix an arbitrary $s$: $0\leq s\leq N-1$. Let us use the induction argument.
For $k=0$ relation~(\ref{f3_17}) holds. Assume that it is true for $k=r\in \mathbb{Z}_+$.
Then
$$ W x^{r+1} \vec e_s = WXW^{-1} Wx^r \vec e_s = \widetilde A x_{rN+s} = x_{(r+1)N+s}. $$
Therefore relation~(\ref{f3_17}) is true.

\noindent
Repeating arguments from~\cite[pp.276-277]{cit_14000_Z} we construct a unitary transformation $V$
which maps $L^2_0(M)$ onto $H$, such that
\begin{equation}
\label{f3_18}
V x^k \vec e_s = x_{kN+s},\qquad 0\leq s\leq N-1;\ k\in \mathbb{Z}_+.
\end{equation}
By~(\ref{f3_17}),(\ref{f3_18}) we conclude that $Wf=Vf$, $f\in L^2_0(M)$.
Therefore $WL^2_0(M) = H$, and $L^2_0(M)=W^{-1}H = L^2(M)$.

\noindent
(ii)$\Leftrightarrow$(iii): This equivalence was established before the statement of the Theorem.

\noindent
(ii)$\Rightarrow$(iv):
Let $M=(m_{k,j})_{k,j=0}^{N-1}$ has form~(\ref{f3_12}) where $\widehat E_t$
is a left-continuous resolution of unity of a self-adjoint operator $\widehat A\supseteq A$ in $H$.
Then
\begin{equation}
\label{f3_19}
( R_\lambda (\widehat A) x_{Nk+r}, x_{Nl+s} )_H =
( R_\lambda (\widehat A) \widehat A^k x_{r}, \widehat A^l x_{s} )_H =
( \widehat A^{k+l} R_\lambda (\widehat A) x_{r}, x_{s} )_H
\end{equation}
\begin{equation}
\label{f3_20}
= \int_\mathbb{R} \frac{ t^{k+l} }{ t-\lambda } d( \widehat E x_r, x_s )_H
= \int_\mathbb{R} \frac{ t^{k+l} }{ t-\lambda } dm_{r,s},\qquad 0\leq r,s\leq N-1;\ k,l\in \mathbb{Z}_+.
\end{equation}
Therefore for $D_\lambda := R_\lambda (\widehat A)$ condition~(iv) holds.

\noindent
(iv)$\Rightarrow$(ii):
Let $E_{U,t}$ be the left-continuous orthogonal resolution of unity of $A_U$.
Observe that $D_\lambda$ is the resolvent function of the self-adjoint operator $A_U\supseteq A$
in $H$.
Using~(\ref{f3_10}) we may write
$$ \int_\mathbb{R} \frac{ 1 }{ x-\lambda } d(E_{U,t}x_r,x_s)_H =
( R_\lambda (A_U) x_{r}, x_{s} )_H =
( D_\lambda x_{r}, x_{s} )_H $$
\begin{equation}
\label{f3_22}
= \int_\mathbb{R} \frac{ 1 }{ x-\lambda } dm_{r,s},\
0\leq r,s\leq N-1.
\end{equation}
Therefore $M= ((E_{U,t}x_r,x_s)_H)_{r,s=0}^{N-1}$. Hence, $M$ is a canonical solution of the moment problem.
$\Box$

\section{Density of polynomials: the case (B).}

Let $\sigma$ be a
non-negative measure on $\mathfrak{B}(\Pi)$, such that
\begin{equation}
\label{f4_1}
\int_\Pi x^m d\sigma < \infty,\quad m\in \mathbb{Z}_+.
\end{equation}
Set
\begin{equation}
\label{f4_2}
s_{m,n} := \int_\Pi x^m e^{in\varphi} d\sigma,\qquad m\in \mathbb{Z_+},\ n\in \mathbb{Z},
\end{equation}
and consider the Devinatz moment problem with moments $\{ s_{m,n} \}_{m\in \mathbb{Z}_+, n\in \mathbb{Z}}$.
Since the moment problem has a solution, for arbitrary
complex numbers $\alpha_{m,n}$ (where all but finite numbers are zeros) we have~\cite{cit_1150_Z}
\begin{equation}
\label{f4_3}
\sum_{m,k=0}^\infty \sum_{n,l=-\infty}^\infty \alpha_{m,n}\overline{\alpha_{k,l}} s_{m+k,n-l} \geq 0.
\end{equation}
There exists a Hilbert space $H$ and a sequence $\{ x_{m,n} \}_{m\in \mathbb{Z}_+, n\in \mathbb{Z}}$ in $H$,
such that $\mathop{\rm span}\nolimits \{ x_{m,n} \}_{m\in \mathbb{Z}_+, n\in \mathbb{Z}} = H$,
and~\cite{cit_1150_Z}
\begin{equation}
\label{f4_4}
(x_{m,n},x_{k,l})_H = s_{m+k,n-l},\qquad m,k\in \mathbb{Z}_+,\ n,l\in \mathbb{Z}.
\end{equation}
Let $A_0$, $B_0$ be linear operators and $J_0$ be an antilinear operator, with
$D(A_0) = D(B_0) = D(J_0) = \mathop{\rm Lin}\nolimits \{ x_{m,n} \}_{m\in \mathbb{Z}_+, n\in \mathbb{Z}}$, defined by
equalities
$$ A_0 x_{m,n} =  x_{m+1,n},\ B_0 x_{m,n} = x_{m,n+1},\ J_0 x_{m,n} = x_{m,-n},\qquad
m\in \mathbb{Z}_+,\ n\in \mathbb{Z}.  $$
In~\cite{cit_1150_Z} it was shown that these operators are correctly defined,
$A_0$ is symmetric and $B_0$ is isometric. Operators $A=\overline{A_0}$
and $B=\overline{B_0}$ are commuting closed symmetric and unitary operators, respectively.
The operator $J_0$ extends by continuity to a conjugation $J$ in $H$.

In~\cite{cit_1150_Z} it was proved that an arbitrary solution $\mu$ of the Devinatz moment problem
has the following form:
\begin{equation}
\label{f4_5}
\mu (\delta)= ((\mathbf{E}\times F)(\delta) x_{0,0}, x_{0,0})_H,\qquad \delta\in \mathfrak{B}(\Pi),
\end{equation}
where $F$ is the spectral measure of $B$, $\mathbf{E}$ is a spectral measure of $A$ which commutes with
$F$. By $((\mathbf{E}\times F)(\delta) x_{0,0}, x_{0,0})_H$ we mean the non-negative Borel measure on
$\Pi$ which is obtained by the Lebesgue continuation procedure from the following
non-negative measure on rectangles
\begin{equation}
\label{f4_6}
((\mathbf{E}\times F)(I_x\times I_\varphi) x_{0,0}, x_{0,0})_H :=
( \mathbf{E}(I_x) F(I_\varphi) x_{0,0}, x_{0,0})_H,
\end{equation}
where $I_x\subset \mathbb{R}$, $I_\varphi\subseteq [-\pi,\pi)$ are arbitrary intervals.

\noindent
On the other hand, for an arbitrary spectral measure $\mathbf{E}$ of $A$ which commutes with the
spectral measure $F$ of $B$, by relation~(\ref{f4_5}) there corresponds a solution of the Devinatz moment
problem.
The correspondence between the spectral measures of $A$ which commute with the spectral measure of
$B$ and solutions of the Devinatz moment problem is bijective.

Recall the following definition~\cite{cit_1150_Z}:
\begin{dfn}
\label{d4_1}
A solution $\mu$ of the Devinatz moment problem~(\ref{f1_1}) is said to be {\bf canonical}
if it is generated by relation~(\ref{f4_5}) where $\mathbf{E}$ is an {\bf orthogonal}
spectral measure of $A$ which commutes with the spectral measure of $B$. Orthogonal spectral measures
are those measures which are the spectral measures of self-adjoint extensions of $A$ inside $H$.
\end{dfn}
We also need some objects introduced in~\cite{cit_1150_Z} to formulate a description of all canonical
solutions. Set $V_A := (A+iE_H)(A-iE_H)^{-1}$, and
\begin{equation}
\label{f4_7}
H_1 := \Delta_A(i),\ H_2 := H\ominus H_1,\ H_3:= \Delta_A(-i),\ H_4 := H\ominus H_3.
\end{equation}
The restriction $B_{H_2}$ of $B$ to $H_2$ is unitary, and by the
Godi\v{c}-Lucenko Theorem it has a
representation: $B_{H_2} = KL$,
where $K$ and $L$ are some conjugations in $H_2$.
Set $U_{2,4} := JK$.
Let $F_2=F_2(\delta)$, $\delta\in \mathfrak{B}([-\pi,\pi))$, be the spectral
measure of the operator $B_{H_2}$ in $H_2$.
Let $\mu$ be a scalar non-negative measure with a type which coincides with
the spectral type of the measure $F_2$.
Let $N_2$ be the multiplicity function of the measure $F_2$. Then there exists a unitary transformation $W$
of the space $H_2$ on the direct integral $\mathcal{H}=\mathcal{H}_{\mu,N_2}$ such that
\begin{equation}
\label{f4_8}
W B_{H_2} W^{-1} = Q_{e^{iy}},
\end{equation}
where $Q_{e^{iy}}: g(y) \mapsto e^{iy} g(y)$.
Denote by $\mathbf{D}(B;H_2)$ a set of all unitary decomposable operators in $\mathcal{H}$.

In relation~(\ref{f4_5}), canonical solutions correspond to those spectral measures $\mathbf{E}$
which are spectral measures of self-adjoint operators $\widehat A$ of the following form:
\begin{equation}
\label{f4_9}
\widehat A = iE_H + 2( V_A\oplus U_{2,4} W^{-1} V_2 W - E_H )^{-1},
\end{equation}
where $V_2\in \mathbf{D}(B;H_2)$. The correspondence between all operators $V_2\in \mathbf{D}(B;H_2)$
and all canonical solutions is bijective~\cite{cit_1150_Z}.

\begin{thm}
\label{t4_1}
Let $\sigma$ be a non-negative measure on $\mathfrak{B}(\Pi)$, such that
relation~(\ref{f4_1}) holds. Let $L^2_0(\sigma)$ be the closure in $L^2(\sigma)$ of a set of all
power-trigonometric polynomials~(\ref{f1_3}).
Consider the Devinatz moment problem with moments $\{ s_{m,n} \}_{m\in \mathbb{Z}_+, n\in \mathbb{Z}}$
defined by~(\ref{f4_2}). Consider
a Hilbert space $H$ and a sequence $\{ x_{m,n} \}_{m\in \mathbb{Z}_+, n\in \mathbb{Z}}$ in $H$,
such that $\mathop{\rm span}\nolimits \{ x_{m,n} \}_{m\in \mathbb{Z}_+, n\in \mathbb{Z}} = H$, and
relation~(\ref{f4_4}) holds.
The following conditions are equivalent:
\begin{itemize}
\item[(i)] $L^2_0(\sigma) = L^2(\sigma)$;

\item[(ii)] $\sigma$ is a canonical solution of the Devinatz moment problem;

\item[(iii)] $\sigma$ is generated by relation~(\ref{f4_5}), where $\mathbf{E}$ is the
spectral function of $\widehat A$ which has the form~(\ref{f4_9})
with an operator $V_2\in \mathbf{D}(B;H_2)$.

\item[(iv)] For every $\lambda\in \mathbb{C}_+$, there exists a linear bounded operator $D_\lambda$
in $H$ such that
\begin{equation}
\label{f4_10}
( D_\lambda x_{m,n}, x_{m',n'} )_H = \int_\Pi \frac{ x^{m+m'} e^{i(n-n')\varphi} }{ x-\lambda } d\sigma,\quad
m,m'\in \mathbb{Z}_+,\ n,n'\in \mathbb{Z},
\end{equation}
which is invertible, and
\begin{equation}
\label{f4_11}
( (E_H + 2i D_i)^k x_{0,n}, x_{0,0} )_H =
\int_\Pi \left( \frac{x+i}{x-i} \right)^k e^{in\varphi} d\sigma,\quad
n,k\in \mathbb{Z};
\end{equation}
\begin{equation}
\label{f4_12}
D_\lambda^{-1} + \lambda E_H \equiv \widehat A,
\end{equation}
where $\widehat A$ has the form~(\ref{f4_9}) with an operator $V_2\in \mathbf{D}(B;H_2)$.
\end{itemize}
\end{thm}
{\bf Proof. }
(i)$\Rightarrow$(ii):
This implication was proved in~\cite{cit_1150_Z} (see considerations before References).

\noindent
(ii)$\Rightarrow$(i):
Let $\sigma$ has form~(\ref{f4_5}), where $\mathbf{E}$
is the spectral function a self-adjoint operator $\widehat A\supseteq A$ in $H$, which commutes with $B$.
Since $\widehat A x_{m,n} = A x_{m,n} = x_{m+1,n}$, $m\in \mathbb{Z}_+$, $n\in \mathbb{Z}$,
by an induction argument we get
\begin{equation}
\label{f4_13}
\widehat A^r x_{m,n} = x_{m+r,n},\qquad m,r\in \mathbb{Z}_+,\ n\in \mathbb{Z}.
\end{equation}
Therefore
$$ \widehat A^r B^l x_{0,0} = \widehat A^r x_{0,l} = x_{r,l},\quad r,l\in \mathbb{Z}_+. $$
We conclude that
\begin{equation}
\label{f4_14}
\mathop{\rm span}\nolimits \{ \widehat A^m B^n x_{0,0},\ m\in \mathbb{Z}_+,\ n\in \mathbb{Z} \} = H.
\end{equation}
Thus, $(\widehat A,B)$ has a spectrum of multiplicity $1$.
By Theorem~\ref{t2_1} there exists a unitary transformation $W$ which maps $L^2(\sigma)$ onto $H$ such that:
\begin{equation}
\label{f4_15}
W^{-1} \widehat A W = X,\ W^{-1} B W = U
\end{equation}
\begin{equation}
\label{f4_16}
W 1 = x_{0,0},
\end{equation}
where $X:\ f(x,\varphi)\mapsto xf(x,\varphi)$ and
$U:\ f(x,\varphi)\mapsto e^{i\varphi} f(x,\varphi)$
in $L^2(\sigma)$.
Let us check that
\begin{equation}
\label{f4_17}
W x^m = x_{m,0},\qquad m\in \mathbb{Z}_+.
\end{equation}
For $m=0$ it is true. Assume that it is true for $r\in \mathbb{Z}_+$. Then
$$ W x^{r+1} = WXW^{-1}Wx^r = \widehat A x_{r,0} = x_{r+1,0}, $$
and therefore~(\ref{f4_17}) holds.
Let us prove that
\begin{equation}
\label{f4_18}
W x^m e^{in\varphi} = x_{m,n},\qquad m\in \mathbb{Z}_+,\ n\in \mathbb{Z}.
\end{equation}
Fix an arbitrary $m\in \mathbb{Z}_+$.
For $n=0$ relation~(\ref{f4_18}) holds. Assume that it is true for $n=r\in \mathbb{Z}_+$.
Then
$$ W e^{ i(r+1)\varphi } x^m = WUW^{-1} W e^{ir\varphi} x^m = B x_{m,r} = x_{m,r+1}. $$
On the other hand, assume that~(\ref{f4_18}) holds for $n=-r$, $r\in \mathbb{Z}_+$. Then
$$  W e^{ i(-r-1)\varphi } x^m = WU^{-1}W^{-1} W e^{-ir\varphi} x^m = B^{-1} x_{m,-r} = x_{m,-r-1}. $$
Therefore relation~(\ref{f4_18}) is true.

\noindent
Repeating arguments from the beginning of the Proof of Theorem 3.1 in~\cite{cit_1150_Z}
we construct a unitary transformation $V$
which maps $L^2_0(\sigma)$ onto $H$, such that
\begin{equation}
\label{f4_19}
V x^m e^{in\varphi} = x_{m,n},\qquad m\in \mathbb{Z}_+,\ n\in \mathbb{Z}.
\end{equation}
By~(\ref{f4_18}),(\ref{f4_19}) we conclude that $Wf=Vf$, $f\in L^2_0(\sigma)$.
Therefore $WL^2_0(\sigma) = H$, and $L^2_0(\sigma)=W^{-1}H = L^2(\sigma)$.

\noindent
(ii)$\Leftrightarrow$(iii): This equivalence was established in~\cite[Theorem 3.2]{cit_1150_Z} and
discussed before the statement of the Theorem.

\noindent
(ii)$\Rightarrow$(iv):
Let $\sigma$ has form~(\ref{f4_5}), where $\mathbf{E}$
is the spectral function a self-adjoint operator $\widehat A\supseteq A$ in $H$, which commutes with $B$.
By~considerations before the statement of the Theorem we obtain that $\widehat A$ has
the form~(\ref{f4_9}) with an operator $V_2\in \mathbf{D}(B;H_2)$.
Then
$$ \left( R_\lambda(\widehat A) x_{m,n}, x_{m',n'} \right)_H
= \left( R_\lambda (\widehat A) \widehat A^m B^n x_{0,0}, \widehat A^{m'} B^{n'} x_{0,0} \right)_H $$
$$ = \left( B^{n-n'} \widehat A^{m+m'} R_\lambda (\widehat A)  x_{0,0}, x_{0,0}
\right)_H
= \int_\Pi
\frac{ x^{m+m'} e^{i(n-n')\varphi} }{ t-\lambda } d( \mathbf{E}\times F) x_{0,0}, x_{0,0} )_H $$
\begin{equation}
\label{f4_21}
= \int_\Pi
\frac{ x^{m+m'} e^{i(n-n')\varphi} }{ t-\lambda } d\sigma,\qquad m,m'\in \mathbb{Z}_+,\
n,n'\in \mathbb{Z};
\end{equation}
$$ \left( (E_H + 2i R_i(\widehat A))^k x_{0,n}, x_{0,0} \right)_H
= \left( (E_H + 2i R_i(\widehat A))^k B^n x_{0,0}, x_{0,0} \right)_H $$
$$ = \int_\Pi \left( \frac{x+i}{x-i} \right)^k
e^{in\varphi} d( \mathbf{E}\times F) x_{0,0}, x_{0,0} )_H $$
\begin{equation}
\label{f4_21_1}
= \int_\Pi \left( \frac{x+i}{x-i} \right)^k
e^{ in\varphi } d\sigma,\qquad k,n\in \mathbb{Z}.
\end{equation}
Therefore for $D_\lambda := R_\lambda (\widehat A)$ condition~(iv) holds.

\noindent
(iv)$\Rightarrow$(ii):
Observe that $D_\lambda$ is the resolvent function of the self-adjoint operator $\widehat A\supseteq A$
in $H$ which commutes with $B$. Let $\widehat E$ be the spectral function of $\widehat A$.
Using~(\ref{f4_11}) we may write
$$ \int_\Pi
\left( \frac{x+i}{x-i} \right)^k
e^{ in\varphi }
d( (\widehat E\times F) x_{0,0},x_{0,0})_H
= ( (E_H + 2i R_i(\widehat A))^k B^n x_{0,0}, x_{0,0} )_H = $$
\begin{equation}
\label{f4_22}
= ( D_\lambda (E_H + 2i D_i)^k x_{0,n}, x_{0,0} )_H
= \int_\Pi \left( \frac{x+i}{x-i} \right)^k e^{ in\varphi } d\sigma,\quad
n,k\in \mathbb{Z}.
\end{equation}
Repeating arguments from the Proof of Theorem~3.1~\cite{cit_1150_Z}, we easily obtain that
$\sigma= ( (\widehat E\times F) x_{0,0},x_{0,0} )_H$.
Hence, $\sigma$ is a canonical solution of the moment problem.
$\Box$

\begin{center}
{\large\bf On the density of polynomials in some $L^2(M)$ spaces.}
\end{center}
\begin{center}
{\bf S.M. Zagorodnyuk}
\end{center}

In this paper we study the density of polynomials in some $L^2(M)$ spaces. Two choices of the measure $M$ and
polynomials are considered:
1) a $(N\times N)$ matrix non-negative Borel measure on $\mathbb{R}$ and vector-valued polynomials
$p(x) = (p_0(x),p_1(x),...,p_{N-1}(x))$, $p_j(x)$ are complex polynomials, $N\in \mathbb{N}$;
2) a scalar non-negative Borel measure in a strip
$\Pi = \{ (x,\varphi):\ x\in \mathbb{R}, \varphi\in [-\pi,\pi) \} $, and
power-trigonometric polynomials:
$p(x,\varphi) = \sum_{m=0}^\infty \sum_{n=-\infty}^\infty \alpha_{m,n} x^m e^{in\varphi}$,
$\alpha_{m,n}\in \mathbb{C}$, where all but finite number of $\alpha_{m,n}$ are zeros.
We prove that polynomials are dense in $L^2(M)$ if and only if $M$ is a canonical solution of the
corresponding moment problem.
Using descriptions of canonical solutions, we get conditions for the density of polynomials in $L^2(M)$.
For this purpose, we derive a model for commuting self-adjoint and unitary operators with
a spectrum of a finite multiplicity.

}
\end{document}